\documentclass[12pt]{article}
\usepackage{amssymb}
\usepackage{amsmath}
\usepackage{indentfirst,latexsym}
\usepackage{mathrsfs}
\usepackage{graphicx}
\usepackage{fancyhdr}

\date{}
\begin{document}

\title{The signed Eulerian numbers on involutions}
\author{Marilena Barnabei, Flavio Bonetti, and Matteo Silimbani \thanks{
Dipartimento di Matematica - Universit\`a di Bologna}} \maketitle

\noindent {\bf Abstract.} We define an analogue of signed Eulerian
numbers $f_{n,k}$ for involutions of the symmetric group and
derive some combinatorial properties of this sequence. In
particular, we exhibit both an explicit formula and a recurrence
for $f_{n,k}$ arising from the properties of its generating
function.\newline



\section{Introduction}

\noindent Let $\sigma$ be a permutation in $S_n$. We say that
$\sigma$ has a \emph{descent} at position $i$ whenever
$\sigma(i)>\sigma(i+1)$. Analogously, we say that $\sigma$ has a
\emph{rise} at position $i$ whenever $\sigma(i)<\sigma(i+1)$. The
number of descents (respectively rises) of a permutation $\sigma$
is denoted by $des(\sigma)$ (resp. $ris(\sigma)$). The polynomial
$$A_n(t)=\sum_{\sigma\in S_n} t^{des(\sigma)}=\sum_{k=0}^{n-1}a_{n,k}t^k,$$ is
known as the \emph{Eulerian polynomial}, and the integers
$a_{n,k}$, i.e., the number of permutations $\sigma\in S_n$ with
$des(\sigma)=k$, are called the \emph{Eulerian numbers}.

\noindent We recall that the Eulerian numbers satisfy the property
$a_{n,k}=a_{n,n-1-k}$, that implies
\begin{equation}\sum_{\sigma\in S_n}
t^{des(\sigma)}=\sum_{\sigma\in S_n}
t^{ris(\sigma)}\label{citare}.\end{equation} The study of the
distribution of the descent statistic has been carried out both in
the case of the symmetric group and of some particular subsets of
permutations. For example, Eulerian distribution on the set of
involutions $\mathscr{I}_n\subseteq S_n$ has been deeply
investigated by several authors (\cite{gess}, \cite{df1},
\cite{zeng}, and \cite{bbsy}).\newline

\noindent Loday \cite{lod}, in his study of the cyclic homology of
commutative algebras, introduced the sequence $(b_{n,k})$ of
\emph{signed Eulerian numbers}, namely, the coefficients of the
polynomial
$$B_n(t)=\sum_{\sigma\in S_n} \textrm{sgn}(\sigma)t^{ris(\sigma)}=\sum_{k=0}^{n-1} b_{n,k}t^k.$$
In \cite{df}, D\'esarm\'enien and Foata showed several
combinatorial properties satisfied by the integers $b_{n,k}$ by
exploiting their relations with the Eulerian numbers. More
recently, Tanimoto (see, e.g., \cite{tan}) developedthe study of
signed Eulerian numbers from a different point of view.\newline

\noindent In this paper we study the signed Eulerian numbers on
 involutions, i.e., the coefficients of the \emph{signed
Eulerian polynomial on involutions}
$$F_n(t)=\sum_{\sigma\in\mathscr{I}_n} \textrm{sgn}\sigma t^{ris(\sigma)}=\sum_{k=0}^{n-1}
f_{n,k}t^k.$$
\newline

\noindent To this aim, we exploit a map introduced in \cite{bbsw}
that associates an involution with a family of generalized
involutions, which share the shape of the corresponding
semistandard Young tableau. This map yields an explicit formula
for the number of involutions in $\mathscr{I}_n$ with a given
number of rises. The above map yields also an analogue of
Worpitzky Identity for the sequence $f_{n,k}$:
$$\sum_{j=0}^{s-1}{n+j\choose j}f_{n,s-j-1}=\sum_{j=0}^{\left\lfloor\frac{n}{2}\right\rfloor} (-1)^j{{s+1\choose 2}+j-1\choose j}{s\choose n-2j}$$
and, hence, an explicit formula for the signed Eulerian numbers on
involutions.

\noindent This formula gives the following expression for the
generating function of the signed Eulerian polynomials on
involutions:
$$\sum_{n\geq 0} F_n(t)\frac{u^n}{(1-t)^{n+1}}=\sum_{r\geq 0} t^r\frac{(1+u)^{r+1}}{(1+u^2)^{r+2\choose 2}}.$$
Applying the Maple package $ZeilbergerRecurrence(T,n,k,s,0..n)$ to
this last identity, we find a recurrence satisfied by the integers
$f_{n,k}$.

\noindent Moreover, we exhibit an explicit formula, a recurrence,
and a generating function for the sequence $F_n(1)$, with
$n\in\mathbb{N}$, namely, the difference between the number of
even and odd involutions in $\mathscr{I}_n$.

\section{The signed Eulerian numbers}

\noindent In order to investigate the combinatorial properties of
the signed Eulerian polynomial on involutions, we exploit the
relations between involutions and particular biwords, called
generalized involutions, introduced in \cite{bbsy}.\newline

\noindent A \emph{generalized involution} of length $n$ is defined
to be a biword:
$$\alpha={x\choose y}=\left(\begin{array}{cccc}
x_1&x_2&\cdots&x_n\\
y_1&y_2&\cdots&y_n\\
\end{array}\right),$$
such that: \begin{itemize}\item for every $1\leq i\leq n$, there
exists an index $j$ with $x_i=y_j$ and $y_i=x_j$, \item $x_i\leq
x_{i+1}$,
\item $x_i=x_{i+1}\Longrightarrow y_i\geq
y_{i+1}$.\end{itemize}

\noindent We say that an integer $a$ is a \emph{repetition} of
multiplicity $r$ for the generalized involution $\alpha$ if
$$x_i=y_i=x_{i+1}=y_{i+1}=\cdots=x_{i+r-1}=y_{i+r-1}=a.$$
\newline

\noindent In \cite{bbsw}, the authors introduced a map $\Pi$ from
the set of generalized involutions of length $n$ to the set of
involutions $\mathscr{I}_n$ defined as follows: if

$$\alpha=\left(\begin{array}{cccc}
x_1&x_2&\cdots&x_n\\
y_1&y_2&\cdots&y_n\\
\end{array}\right),$$
then $\Pi(\alpha)$ is the involution $\sigma$
$$\sigma=\left(\begin{array}{cccc}
1&2&\cdots&n\\
y'_1&y'_2&\cdots&y'_n\\
\end{array}\right),$$
where $y'_i=1$ if $y_i$ is the least symbol occurring in the word
$y$, $y'_j=2$ if $y_j$ is the second least symbol in $y$ and so
on. In the case $y_i=y_j$, with $i>j$, we consider $y_i$ to be
less then $y_j$. We call the involution $\sigma=\Pi(\alpha)$ the
\emph{polarization} of $\alpha$.
\newline

\noindent We denote by Gen$_m(\sigma)$ the set of generalized
involutions of length $n$, with symbols taken from $[m]$, whose
polarization is $\sigma$. Then, we have the following result
proved in \cite{bbsy}:

\newtheorem{yama}{Proposition}
\begin{yama}\label{nlac}
Let $\sigma\in \mathscr{I}_n$ be an involution with $t$ rises.
Then,
\begin{equation}|\textrm{Gen}_m(\sigma)|={n+m-t-1\choose n}.
\label{consalite}\end{equation}
\end{yama}
\begin{flushright}
$\diamond$
\end{flushright}

\noindent We recall that the sign of an involution $\sigma\in
\mathscr{I}_n$ is determined by the number $fix(\sigma)$ of fixed
points of $\sigma$. More precisely:
$$\textrm{sgn}(\sigma)=(-1)^{\frac{n-fix(\sigma)}{2}}.$$
It is well known \cite{schutz} that the integer $fix(\sigma)$
depends only on the shape of the tableau $T_{\sigma}$. On the
other hand, it is immediately seen that the semistandard tableau
associated with any generalized involution in Gen$_m(\sigma)$ by
the Robinson-Schensted-Knuth algorithm has the same shape as the
standard tableau $T_{\sigma}$ associated with $\sigma$. We are
interested in counting separately the number of even and odd
involutions with a given number of rises.  Consequently, we can
define an \emph{even generalized involution} to be a generalized
involution whose polarization is an even involution.

\noindent Remark that the sign of a generalized involution depends
only on the multiplicity of its repetitions. In fact, given a
generalized involution $\alpha$, we define $gfix(\alpha)$ to be
the number of repetitions of odd multiplicity of $\alpha$ (note
that, if $\alpha\in \mathscr{I}_n$, then
$gfix(\alpha)=fix(\alpha)$). It is easy to verify that a
generalized involution $\alpha$ of length $n$ is even if and only
if
$$\frac{n-gfix(\sigma)}{2}$$
is even. This remark allows to compute explicitly the number
$a^+_{n,m}$ of generalized involutions of length $n$ over $[m]$:

\newtheorem{solopari}[yama]{Proposition}
\begin{solopari}\label{mucho}
We have:
$$a^+_{n,m}=\sum_{h=0}^{\left\lfloor\frac{n}{4}\right\rfloor} {m\choose n-4h}{{m+1\choose 2}+2h-1\choose 2h}.$$
\end{solopari}

\noindent \emph{Proof} Let $\alpha$ be a generalized involution.
We remarked that $\alpha$ is even whenever
$\frac{n-gfix(\sigma)}{2}$ is even. We want to count generalized
involutions with $\frac{n-gfix(\sigma)}{2}=2h$ with fixed $h$. We
can choose the $n-4h$ fixed points of $\alpha$ in ${m\choose
n-4h}$ ways. Then, we choose $2h$ biletters $(i,j)$ in the set
$B=\{(i,j)|1\leq i\leq j\leq m\},$ whose cardinality is
${m+1\choose 2}$. We complete the generalized involution by adding
a biletter $(j,i)$ for every biletter $(i,j)$ previously chosen,
hence getting an even generalized involution.
\begin{flushright}
$\diamond$
\end{flushright}

\noindent On the other hand, recalling that the total number of
generalized involutions of length $n$ over $[m]$ is (see, e.g.,
\cite{bbsy})
$$a_{n,m}=\sum_{h=0}^{\left\lfloor\frac{n}{2}\right\rfloor} {m+n-2h-1\choose n-2h}{{m\choose 2}+h-1\choose h},$$
we get an explicit formula for the number $a^-_{n,m}$ of odd
generalized involutions of length $n$ over $[m]$:
$$a^-_{n,m}=\sum_{h=0}^{\left\lfloor\frac{n}{2}\right\rfloor} {m+n-2h-1\choose n-2h}{{m\choose 2}+h-1\choose h}-\sum_{h=0}^{\left\lfloor\frac{n}{4}\right\rfloor} {m\choose n-4h}{{m+1\choose 2}+2h-1\choose 2h}.$$
We remark that, given $\sigma\in\mathscr{I}_n$, all generalized
involutions in Gen$_m(\sigma)$ have the same sign as $\sigma$.
Hence, Proposition \ref{nlac} implies that
\begin{equation}a^+_{n,m}=\sum_{k=0}^{m-1}{n+k\choose k} f^+_{n,m-k-1}\label{piu}\end{equation}
\begin{equation}a^-_{n,m}=\sum_{k=0}^{m-1}{n+k\choose k} f^-_{n,m-k-1}\label{meno},\end{equation}
where $f^+_{n,k}$ and $f^-_{n,k}$ denote the number of of positive
and negative involutions in $\mathscr{I}_n$ with $k$ rises,
respectively. These identities allow to state the following:

\newtheorem{segnati}[yama]{Theorem}
\begin{segnati}\label{hnc}
We have: \begin{equation}f_{n,k}=\sum_{m=0}^{k+1}(-1)^{k-m+1}{n+1
\choose
k-m+1}\sum_{j=0}^{\left\lfloor\frac{n}{2}\right\rfloor}(-1)^j{{m+1\choose
2}+j-1\choose j}{m\choose n-2j}.\label{theo}\end{equation}
\end{segnati}

\noindent \emph{Proof} Set $\hat{a}_{n,m}=a^+_{n,m}-a^-_{n,m}$.
Then, Identities (\ref{piu}) and (\ref{meno}) imply
\begin{equation}\hat{a}_{m,n}=\sum_{k=0}^{m-1}{n+k\choose k} f_{n,m-k-1}.\label{saraworp}\end{equation}
By inversion, we have:
$$f_{n,k}=\sum_{m=0}^{k+1}(-1)^{k-m+1}{n+1 \choose k-m+1}\hat{a}_{n,m}=$$
$$=\sum_{m=0}^{k+1}(-1)^{k-m+1}{n+1 \choose
k-m+1}\left(2\sum_{j=0}^{\left\lfloor\frac{n}{4}\right\rfloor}
{m\choose n-4j}{{m+1\choose 2}+2j-1\choose
2j}\right.$$$$\left.-\sum_{j=0}^{\left\lfloor\frac{n}{2}\right\rfloor}
{m+n-2j-1\choose n-2j}{{m\choose 2}+j-1\choose j}\right)$$ that is
equivalent to (\ref{theo}).
\begin{flushright}
$\diamond$
\end{flushright}

\noindent Remark that Identity (\ref{saraworp}) gives the
following analogue of the Worpitzky Identity:
$$\sum_{j=0}^{s-1}{n+j\choose j}f_{n,s-j-1}=\sum_{j=0}^{\left\lfloor\frac{n}{2}\right\rfloor} (-1)^j{{s+1\choose 2}+j-1\choose j}{s\choose n-2j}.$$

\section{The signed Eulerian polynomial on involutions}

\noindent We define the $n$-th \emph{signed Eulerian polynomial
for involutions} to be the polynomial:
$$F_n(t)=\sum_{\sigma\in\mathscr{I}_n} \textrm{sgn}(\sigma) t^{ris(\sigma)}=\sum_{k=0}^{n-1}
f_{n,k}t^k.
$$
As shown in the previous section, the relation between involutions
and generalized involutions is crucial in our analysis. We denote
by
$$R_n(t)=\sum_{m\geq 0} \hat{a}_{n,m}t^m$$
$$C_m(t)=\sum_{n\ge 0} \hat{a}_{n,m}t^n$$
the row and column generating functions of the array
$\hat{a}_{n,m}$. As seen in the proof of Theorem \ref{hnc}, we
have
$$\hat{a}_{n,m}=\sum_{j=0}^{\left\lfloor\frac{n}{2}\right\rfloor}(-1)^j{{m+1\choose 2}+j-1\choose j}{m\choose n-2j}.$$
and hence
$$C_m(t)=\frac{(1+u)^m}{(1+u^2)^{m+1\choose 2}}.$$

\newtheorem{polygen}[yama]{Theorem}
\begin{polygen}
The polynomial $F_n(t)$ satisfies the identity
\begin{equation} \sum_{n\geq 0}F_n(t)\frac{u^n}{(1-t)^{n+1}}=\sum_{m\geq 0}t^m\frac{(1+u)^{m+1}}{(1+u^2)^{m+2\choose 2}}.\label{perguo}\end{equation}
\end{polygen}

\noindent \emph{Proof} The binomial relation between the sequences
$f_{n,k}$ and $\hat{a}_{n,m}$ yields the following identity:
$$\frac{tF_n(t)}{(1-t)^{n+1}}=R_n(t).$$
Hence:
$$\sum_{n\geq 0}F_n(t)\frac{u^n}{(1-t)^{n+1}}=\sum_{n\geq 0}\sum_{m\geq 0} \hat{a}_{n,m+1}t^m u^n=$$$$=\sum_{m\geq 0} C_{m+1}t^m=\sum_{m\geq 0}t^m\frac{(1+u)^{m+1}}{(1+u^2)^{m+2\choose 2}}$$
as desired.
\begin{flushright}
$\diamond$
\end{flushright}

\noindent We apply the Maple package
$ZeilbergerRecurrence(T,n,k,s,0..n)$ to the sequence $f_{n,k}$,
following along the lines of \cite{zeng}, and we deduce by
Identity (\ref{perguo}) the following recurrence formula:

$$n f_{n,k}=(2+k-n)f_{n-1,k}+(2n-k-1)f_{n-1,k-1}-(n+3k+k^2)f_{n-2,k}+$$
$$+(-2+4k+2k^2-2kn)f_{n-2,k-1}+(2-k-k^2+2kn-n^2)f_{n-2,k-2}+$$
$$+(-n-k^2-2k+2)f_{n-3,k}+(-7+4k+3k^2+2n-2kn)f_{n-3,k-1}+$$
$$+(8-2k-3k^2-2n+4kn-n^2)f_{n-3,k-2}+(-3+k^2+n-2kn+n^2)f_{n-3,k-3}.\vspace{.5cm}$$
\noindent In conclusion, we study the combinatorial properties of
the sequence $F_n(1)$, with $n\in \mathbb{N}$. Obviously, $F_n(1)$
is the difference between the number $i^+_n$ of even involutions
on $n$ objects and the number $i^-_{n}$ of odd involutions. First
of all, we have:

\newtheorem{quinci}[yama]{Proposition}
\begin{quinci}
The evaluation of the polynomial $F_n(t)$ at $1$ is
$$F_n(1)=2\sum_{h=0}^{\left\lfloor\frac{n}{4}\right\rfloor}\frac{n!}{(2h)!(n-4h)!2^{2h}}-
\sum_{h=0}^{\left\lfloor\frac{n}{2}\right\rfloor}\frac{n!}{(n-2h)!2^h}.$$
\end{quinci}

\noindent \emph{Proof} Fix an integer $h\leq
\lfloor\frac{n}{4}\rfloor$. We count the number of involutions
whose cycle decomposition consists of $2h$ transpositions. Choose
a word $w=w_1\cdots w_n$ consisting of distinct letters taken from
$[n]$. We have $n!$ choices for $w$. This word correspond to a
unique even involution $\tau$ with $n-4h$ fixed points defined by
the following conditions:
$$\tau(w_1)=w_2\quad\ldots\quad\tau(w_{4h-1})=w_{4h};$$
$$\tau(w_{4h+j})=w_{4h+j}$$
with $0<j\leq n-4h$. It is easily checked that the involution
$\tau$ arises from $(n-4h)!(2h)!2^{2h}$ different words $w$. These
considerations imply that
$$i^+_n=\sum_{h=0}^{\left\lfloor\frac{n}{4}\right\rfloor}\frac{n!}{(2h)!(n-4h)!2^{2h}}.$$
On the other hand, it is well known that
$$|\mathscr{I}_n|=\sum_{h=0}^{\left\lfloor\frac{n}{2}\right\rfloor}\frac{n!}{(n-2h)!2^h}.$$
Hence, since $F_n(1)=i^+_n-i^-_n=2i^+_n-|\mathscr{I}_n|$, we get
the assertion.
\begin{flushright}
$\diamond$
\end{flushright}

\newtheorem{invorecu}[yama]{Proposition}
\begin{invorecu}
The sequence $F_n(1)$ satisfies the recurrence
$$F_n(1)=F_{n-1}(1)-(n-1)F_{n-2}(1).$$
\end{invorecu}

\noindent \emph{Proof} Consider an even involution
$\tau\in\mathscr{I}_n$. If $\tau(1)=1$, the restriction of $\tau$
on the set $\{2,\ldots,n\}$ is an even involution on $n-1$
objects. If $\tau(1)=j\neq 1$, the restriction of $\tau$ on
$\{2,\ldots,n\}\setminus \{j\}$ is an odd involution on $n-2$
objects. This implies that:
$$i^+_n=i^+_{n-1}+(n-1)i^-_{n-2}$$
and, analogously,
$$i^-_n=i^-_{n-1}+(n-1)i^+_{n-2}.$$
These identities gives immediately the assertion.
\begin{flushright}
$\diamond$
\end{flushright}

\noindent These results allow to deduce the following expression
for the exponential generating function of the sequence $F_n(1)$:
$$\sum_{n\geq 0}\frac{F_n(1)}{n!}t^n=e^{t-\frac{t^2}{2}}.$$

\end{document}